\documentclass[12pt]{article} 
\usepackage{amsmath,amssymb,eucal} 
\usepackage{bbm}  
\font\tenrm=cmr10

\font\cmssl=cmss10 at 12 pt  
\font\bigss=cmssdc10 scaled 2300

\font\cmsslll=cmss10 at 14 pt

\newcommand{\e}{\epsilon}

\newcommand{\D}{\Delta}  
  
\newcommand{\G}{\Gamma}  
  
\newcommand{\I}{\tau}  
  
\renewcommand{\O}{\Omega}

\newcommand{\bR}{\mathbb{R}}  
\newcommand{\bZ}{\mathbb{Z}}  
  
\newcommand{\bN}{\mathbb{N}}

\newcommand{\id}   {{\mathbbm{1}}}

\renewcommand{\square}{\kern1pt\vbox  
               {\hrule height 0.6pt\hbox{\vrule width 0.6pt\hskip 3pt  
    \vbox{\vskip 6pt}\hskip 3pt\vrule width 0.6pt}\hrule height0.6pt}  
    \kern1pt}

\newcommand{\ra}{\rightarrow}

\newtheorem{Th}{Theorem}  
\newtheorem{Prop}{Proposition}  
\newtheorem{Cor}{Corollary}  
\newtheorem{Lem}{Lemma}  
\newtheorem{Def}{Definition}  
\newtheorem{Ex}{Example}  
\newcommand{\bt}{\begin{Th}\ \ }  
\newcommand{\et}{\end{Th}}  
\newcommand{\bp}{\begin{Prop}\ \ }  
\newcommand{\ep}{\end{Prop}}  
\newcommand{\bc}{\begin{Cor}\ \ }  
\newcommand{\ec}{\end{Cor}}  
\newcommand{\bl}{\begin{Lem}\ \ }  
\newcommand{\el}{\end{Lem}}  
\newcommand{\bd}{\begin{Def}\ \ }  
\newcommand{\ed}{\end{Def}}  
  
\newcommand{\pf}{\noindent{\it Proof:\ \ }}  
\newcommand{\qed}{\hfill\square}  
\newcommand{\n}{\nabla}

\newcommand{\be}{\begin{equation}}  
\newcommand{\ee}{\end{equation}}  
  
\newcommand\re[1]{(\ref{#1})}  
\newcommand{\arr}{\begin{array}{rlll}}  
\newcommand{\ea}{\end{array}}  
\newcommand{\bea}{\begin{eqnarray}}  
\newcommand{\eea}{\end{eqnarray}}  
\newcommand{\bean}{\begin{eqnarray*}}  
\newcommand{\eean}{\end{eqnarray*}}  
  
\catcode`@=11  
\@addtoreset{equation}{section}  
\catcode`@=12

\begin{document}  
\begin{titlepage}  
\vskip 1.5 true cm  
\begin{center}  
{\bigss  Geometric structures on Lie groups\\[.5em] with flat bi-invariant metric}  
\vskip 1.0 true cm   
{\cmsslll    Vicente Cort\'es$^1$ and Lars Sch\"afer$^2$} \\[3pt]  
{\tenrm   Department Mathematik$^1$ \\
Universit\"at Hamburg \\ Bundesstra{\ss}e 55 \\D-20146 Hamburg, Germany  \\
cortes@math.uni-hamburg.de

Institut f\"ur Differentialgeometrie$^2$ \\ Leibniz Universit\"at Hannover \\
Welfengarten 1\\ D-30167 Hannover, Germany \\
 schaefer@math.uni-hannover.de}\\[1em]   
October 23, 2008 
\end{center}  
\vskip 1.0 true cm  
%%%%%%%%%%%%%%%%%%%%%%%%%%%%%%%%%%%%%%%%%%%%%%%%%%%%%%%%  
\baselineskip=18pt  
\begin{abstract}  
\noindent  
Let $L\subset V=\bR^{k,l}$ be a maximally isotropic subspace.
It is shown that any simply connected Lie group with a bi-invariant flat 
pseudo-Riemannian metric of signature $(k,l)$ is 
2-step nilpotent and 
is defined by an element $\eta \in \Lambda^3L\subset \Lambda^3V$. 
If $\eta$ is of type $(3,0)+(0,3)$ with respect to a skew-symmetric 
endomorphism $J$ with $J^2=\e Id$, then the Lie group
${\cal L}(\eta )$ is endowed with a left-invariant nearly K\"ahler 
structure if $\e =-1$ and with a left-invariant nearly para-K\"ahler 
structure if $\e =+1$. This construction 
exhausts all complete simply connected 
flat  nearly (para-)K\"ahler manifolds. If $\eta \neq 0$ has rational 
coefficients with respect to some basis, then ${\cal L}(\eta )$
admits a lattice $\Gamma$, and the quotient 
$\Gamma\setminus {\cal L}(\eta )$ is a compact 
inhomogeneous nearly (para-)K\"ahler manifold. The first non-trivial
example occurs in six dimensions. \\
{\it MSC(2000):} 53C50, 53C15.   \\
{\it Keywords:} Flat Lie-groups, bi-invariant metrics,  nearly para-K\"ahler manifolds, flat almost (para-)Hermitian manifolds, almost (para-)complex structures.
\end{abstract}  
\vfill \hrule width 3.0 cm  
{\small \noindent This work was supported by the SFB 676 of the Deutsche Forschungsgemeinschaft. }  
\end{titlepage}  
\tableofcontents
\section*{Introduction}
A pseudo-Riemannian manifold $(M,g)$ endowed with a skew-symmetric 
almost complex structure 
$J$ is called {\cmssl nearly K\"ahler} if the Levi-Civita 
covariant derivative $DJ$ is skew-symmetric, that is $(D_XJ)X=0$ for all 
$X\in TM$. Nearly K\"ahler manifolds with a positive definite metric 
are by now well studied, see \cite{N} and 
references therein. Replacing the equation $J^2=-Id$ by $J^2=Id$ one arrives
at the definition of {\cmssl nearly para-K\"ahler} manifold, see
\cite{I}. This generalises the notion of a para-K\"ahler (or bi-Lagrangian) 
manifold. Such manifolds occur naturally in super-symmetric field theories 
over Riemannian rather than Lorentzian space-times, see \cite{CMMS}. 
In \cite{I} Ivanov and 
Zamkovoy ask for examples of
Ricci-flat nearly para-K\"ahler manifolds 
in six dimensions with $DJ\neq 0$. In this 
paper we will give a classification of flat nearly para-K\"ahler manifolds.
In particular, we will show that there exists a 
compact six-dimensional such manifold with $DJ\neq 0$. 

It is noteworthy that flat nearly para-K\"ahler manifolds $M$ 
provide also solutions of the 
so-called tt*-equations, see \cite{S} and references therein. 
As a consequence, they give rise to a (para-)pluriharmonic map from $M$ into 
the pseudo-Riemannian symmetric space $SO_0(n,n)/GL(n)$.
  
Let $V$ be a pseudo-Euclidian vector space and $\eta\in \Lambda^3V$. 
Contraction with $\eta$ defines a linear map $\Lambda^2V^* \ra V$.
The image of that map is denoted by $\Sigma_\eta$ and is called
the {\cmssl support} of $\eta$.  In the first section we will show that 
any 3-vector $\eta\in \Lambda^3V$ 
with isotropic support defines a simply connected 2-step nilpotent 
Lie group ${\cal L}(\eta )$ with a flat bi-invariant 
pseudo-Riemannian metric $h$ of the same signature as $V$. We prove that 
this exhausts all simply connected Lie groups with a flat bi-invariant metric,
see Theorem \ref{Thm2}. After completion of our article, Oliver Baues,
has kindly informed us about the paper \cite{W}, which already 
contains a version of that result. 

It is shown that the groups $({\cal L}(\eta ),h)$ admit a 
lattice $\Gamma \subset {\cal L}(\eta )$ if $\eta$ has rational coefficients 
with respect to some basis and that the quotient  $M(\eta, \Gamma ) :=
\Gamma\setminus {\cal L}(\eta )$ is a {\em flat compact homogeneous 
pseudo-Riemannian manifold}, see Corollary \ref{Cor}. 
Compact homogeneous flat pseudo-Riemannian manifolds 
were recently classified in independent work by Baues, see \cite{B}.  
It follows from this classification that the above examples 
exhaust all compact homogeneous flat pseudo-Riemannian manifolds. 

Assume now that $\dim V$ is even and that we fix $J\in \mathfrak{so}(V)$
such that $J^2= -Id$ or $J^2= Id$. We denote the corresponding 
left-invariant endomorphism field on the 
group ${\cal L}(\eta )$ again by $J$. 
Assume that $\eta\in \Lambda^3V$ has isotropic support and satisfies, 
in addition,  
$$ \{ \eta_X,J\} := \eta_XJ + J\eta_X=0 \quad\mbox{for all}\quad X\in V,$$
or, equivalently, that $\eta$ has type $(3,0) + (0,3)$.  
Then $({\cal L}(\eta ),h,J)$ is a flat nearly K\"ahler manifold  
if $J^2=-Id$ and a flat nearly para-K\"ahler manifold if $J^2=Id$. 
This follows from the results of \cite{CS} for the former case and is proven 
in the second section of this paper for the latter case, see Theorem 
\ref{1stThm}.  
Moreover it is shown that any complete simply 
connected flat nearly (para-)K\"ahler manifold is of this form, see 
Corollary \ref{lastCor} and \cite{CS}. 
To sum up, we have shown that any simply connected complete flat nearly 
(para-)K\"ahler manifold is a Lie group ${\cal L}(\eta )$ 
with a left-invariant nearly 
(para-)K\"ahler structure and bi-invariant metric.  
Conversely, it follows from unpublished work of Paul-Andi Nagy and 
the first author that a Lie group with a left-invariant nearly 
(para-)K\"ahler structure and bi-invariant metric is necessarily 
flat and is therefore covered by one of our groups ${\cal L}(\eta )$. 
The proof of this statement uses the unicity of the connection 
with totally skew-symmetric torsion preserving the nearly 
(para-)K\"ahler structure and the Jacobi identity. 

Suppose now that $\Gamma \subset {\cal L}(\eta )$ is a lattice.  Then 
the almost (para-)complex structure $J$ on the group ${\cal L}(\eta )$ 
induces an almost (para-)\linebreak[3] complex structure $J$ 
on the compact manifold 
$M=M(\eta,\Gamma )= \Gamma \setminus {\cal L}(\eta )$. Therefore $(M,h,J)$
is a compact nearly (para-)K\"ahler manifold. However, the 
(para-)complex structure  is not ${\cal L}(\eta )$-invariant, unless 
$\eta=0$.        Moreover, for $\eta\neq 0$, $(M,h,J)$ is an inhomogeneous 
nearly (para-)\linebreak[3] K\"ahler manifold, that is, it 
does not admit any transitive group of automorphisms of the nearly 
(para-)K\"ahler structure.   
Since $J$ is not right-invariant, this follows from the fact that ${\mathrm{Isom}}_0(M)$
is obtained from the 
action of ${\cal L}(\eta )$
by right-multiplication on $M$, 
see Corollary \ref{Cor}.  The first such non-trivial flat compact nearly 
para-K\"ahler nilmanifold $M(\eta )= \Gamma \setminus {\cal L}(\eta )$ is
six-dimensional and 
is obtained from a non-zero element $\eta \in \Lambda^3V^+\cong \bR$, where 
$V^+\subset V=\bR^{3,3}$ is the $+1$-eigenspace of $J$.

\section{A class of flat pseudo-Riemannian Lie groups} \label{simplytransSec}
Let $V=(\bR^n,\langle \cdot , \cdot \rangle )$ be the standard 
pseudo-Euclidian vector space of signature $(k,l)$, $n=k+l.$
Using the (pseudo-Euclidian) scalar product we shall identify 
$V\cong V^*$ and $\Lambda^2V\cong \mathfrak{so}(V)$. These
identifications provide the inclusion
$\Lambda^3V \subset V^*\otimes \mathfrak{so}(V)$. Using it 
we consider  a three-vector $\eta\in \Lambda^3V$ as an 
$\mathfrak{so}(V)$-valued one-form. Further  we denote by
$\eta_X\in \mathfrak{so}(V)$ the evaluation of this one-form on a vector
$X\in V$.  
\noindent
The {\cmssl support} of $\eta \in \Lambda^3V$ is defined by
\be \Sigma_{\eta}:= {\rm span}\{\eta_XY \,|\, X,Y \in V\} \subset V.\ee
\bt \label{2stepnp_Thm} 
Each 
$$\eta \in \mathcal{C}(V):= \{\eta \in \Lambda^3V \,|\Sigma_\eta \mbox{  
(totally) isotropic} \}=\underset{L\subset V}{\bigcup} \Lambda^3 L$$
defines a 2-step nilpotent simply transitive subgroup 
 $\mathcal{L}(\eta) \subset \mathrm{Isom}(V),$ where the union runs over 
all maximal  isotropic subspaces. 
The subgroups $\mathcal{L}(\eta)$, $\mathcal{L}(\eta') 
\subset \mathrm{Isom}(V)$ 
associated to $\eta, \eta' \in {\cal C}(V)$ 
are conjugated if and only if $\eta' = g\cdot \eta$ for some element of 
$g \in O(V).$

\et
\pf 
It is easy to see that any three-vector $\eta \in \Lambda^3V$ satisfies
$\eta \in \Lambda^3\Sigma_\eta$, cf.\ \cite{CS} Lemma 7. 
This implies the equation 
$\mathcal{C}(V)=\underset{L\subset V}{\bigcup} \Lambda^3 L$. 
Let an element $\eta\in \mathcal{C}(V)$ be given. 
One can easily show that $\Sigma_\eta$ is isotropic if and only if 
the endomorphisms $\eta_X \in\mathfrak{so}(V)$ 
satisfy $\eta_X \circ \eta_Y =0$ for all $X, Y\in V$, cf.\ \cite{CS} 
Lemma 6. The 2-step nilpotent group  
$$\mathcal{L}(\eta) := \left. \left\{g_X:=\exp \left(
\begin{array}{cc}
\eta_X & X \\
0 &0
\end{array}
\right)=\left(
\begin{array}{cc}
\id + \eta_{X} &  X \\
0 &1
\end{array}
\right) \,\right| \, X \in V \right\} $$
acts simply transitively on $V\cong V\times \{ 1\} \subset V\times \bR$
by isometries: 
$$ \left(
\begin{array}{cc}
\id + \eta_{X} &  X \\
0 &1
\end{array}
\right) 
\left(
\begin{array}{c} 0 \\ 1\end{array}\right) = \left(
\begin{array}{c} X \\ 1\end{array}\right).$$
Let us check that $\mathcal{L}(\eta)$ is a group: Using $\eta_X \circ
\eta_Y=0$ we obtain
\bean 
g_X \cdot g_Y&=&
\left( \begin{array}{cc}
\id + \eta_{X} &  X \\
0 &1
\end{array}
\right)
\left(\begin{array}{cc}
\id + \eta_{Y} &  Y \\
0 &1
\end{array}
\right)= \left(\begin{array}{cc}
\id + \eta_{X}+\eta_Y+ \eta_X\eta_Y  &  X+Y+\eta_XY \\
0 &1
\end{array}
\right)\\
&=&
\left(\begin{array}{cc}
\id + \eta_{X+Y} &   X+Y+\eta_XY \\
0 &1
\end{array}
\right) = g_{X+Y+\eta_XY}.
\eean
In the last step we used $\eta_{\eta_XY}=0,$ which follows from 
$\langle\eta_{\eta_XY}Z,W \rangle= \langle\eta_ZW, \eta_XY\rangle$ for all
$X,Y,Z,W \in V.$ Next we consider $\eta,\eta' \in \mathcal{C}(V)$, $g\in O(V)$. 
The computation 
$$g\mathcal{L}(\eta)g^{-1}=\left\{ \left. \left(
\begin{array}{cc}
\id + g\eta_{X}g^{-1} &  gX \\
0 &1
\end{array}
\right) \,\right| \, X \in V \right\} = \left\{ \left. \left(
\begin{array}{cc}
\id + g\eta_{g^{-1}Y}g^{-1} &  Y \\
0 &1
\end{array}
\right) \,\right| \, Y \in V \right\}
$$
shows that $g\mathcal{L}(\eta)g^{-1}= \mathcal{L}(\eta')$ if and only if 
$\eta'_X = (g \cdot \eta)_X=g \, \eta_{g^{-1}X} \, g^{-1}$ for all $X\in V$. 
\qed

\noindent 
Let ${\cal L} \subset \mathrm{Isom}(V)$ be a simply transitive group. Pulling back
the scalar product on $V$ by the orbit map 
\be {\cal L}\ni g\mapsto 
g0\in V \label{def_obitmap} \ee yields a left-invariant flat pseudo-Riemannian metric $h$ 
on ${\cal L}$. A pair $({\cal L},h)$ consisting of a Lie group $\cal L$ 
and a flat left-invariant  
pseudo-Riemannian metric $h$ on $\cal L$ is called 
a {\cmssl flat pseudo-Riemannian Lie group}. 
\bt \label{Thm2} \begin{itemize}
\item[(i)] 
The class of flat pseudo-Riemannian Lie groups $({\cal L}(\eta ),h)$ 
defined in Theorem \ref{2stepnp_Thm} exhausts all simply connected  flat 
pseudo-Riemannian Lie groups with bi-invariant metric. 
\item[(ii)] A Lie group with 
bi-invariant metric is flat if and only if it is 2-step nilpotent. 
\end{itemize} 
\et 

\pf (i) The group ${\cal L}(\eta )$ associated to
a three-vector $\eta \in \mathcal{C}(V)$ is diffeomorphic to $\bR^n$ by
the exponential map. We have to show that the flat pseudo-Riemannian
metric $h$ on ${\cal L}(\eta )$ is bi-invariant. The Lie algebra 
of ${\cal L}(\eta )$ is identified with the vector space $V$ endowed
with the Lie bracket 
\be \label{LiebracketEqu} [X,Y] := \eta_XY-\eta_YX= 2\eta_XY,\quad X,Y\in V.\ee
The left-invariant metric $h$ on ${\cal L}(\eta )$ corresponds to the 
scalar product $\langle \cdot ,\cdot \rangle$� on $V$. 
Since $\eta\in \Lambda^3V$, the endomorphisms $\eta_X=\frac{1}{2}ad_X$ 
are skew-symmetric. This shows that $h$ is bi-invariant. 

Conversely, let $(V,[\cdot , \cdot ])$ 
be the Lie algebra of a pseudo-Riemannian Lie group of dimension $n$ 
with bi-invariant metric $h$. We can assume that the bi-invariant metric 
corresponds to the standard scalar product $\langle \cdot ,\cdot \rangle$ 
of signature $(k,l)$ on $V$. Let us denote by $\eta_X \in \mathfrak{so}(V)$, 
$X\in V$, the skew-symmetric endomorphism of $V$ which corresponds to the 
Levi-Civita covariant derivative $D_X$ acting on left-invariant vector 
fields. {}From the bi-invariance and the Koszul formula we obtain 
that $\eta_X=\frac{1}{2}ad_X$ and, hence,  $R(X,Y)= -\frac{1}{4}ad_{[X,Y]}$
for the curvature.  The last formula shows that $h$ is flat if and only if
the Lie group is 2-step nilpotent. This proves (ii). To finish the 
proof of (i)  we have to show that, under this assumption, 
$\eta$ is completely skew-symmetric
and has isotropic support. The complete skew-symmetry follows from
$\eta_X=\frac{1}{2}ad_X$ and the bi-invariance. Similarly, using the 
bi-invariance, we have 
$$4\langle \eta_XY,\eta_ZW\rangle =\langle [X,Y],[Z,W]\rangle =  -\langle Y,[X,[Z,W]]\rangle  =0 ,$$ 
since the Lie algebra is 2-step nilpotent. This shows that
$\Sigma_\eta$ is isotropic. 
\qed  

\bc With the above notations, let 
$L\subset V$ be a maximally isotropic subspace. 
The correspondence $\eta \mapsto {\cal L}(\eta )$ defines a 
bijection between the points of the orbit space $\Lambda^3L/GL(L)$ and
isomorphism classes of pairs $({\cal L}, h)$ consisting of a 
simply connected Lie group  
$\cal L$ endowed with a flat
bi-invariant pseudo-Riemannian metric $h$ of signature $(k,l)$.  
\ec  

\bc Any simply connected Lie group ${\cal L}$ with a flat bi-invariant metric 
$h$ of signature 
$(k,l)$ contains a normal subgroup of dimension $\ge \max (k,l)\ge 
\frac{1}{2}\dim V$ which acts by translations on the pseudo-Riemannian
manifold $({\cal L},h)\cong \bR^{k,l}$.   
\ec 
\pf 
Let $\mathfrak{a}:= ker (X\mapsto \eta _X) \subset V$ be the kernel of 
$\eta$. Then $\mathfrak{a} = \Sigma_\eta^\perp$ is co-isotropic and defines 
an Abelian ideal $\mathfrak{a}\subset \mathfrak{l}:=Lie\, {\cal L}\cong V
\cong \bR^{k,l}$. The corresponding normal subgroup 
$A\subset {\cal L}={\cal L}(\eta )$
is precisely the subgroup of translations. So we have shown that  
$\dim A \ge  \max (k,l)\ge \frac{1}{2}\dim V$.  
\qed

\noindent 
{\bf Remarks} 1) The number $\dim \Sigma_\eta$ is an isomorphism invariant of 
the groups ${\cal L}={\cal L}(\eta )$, which is independent of the metric. 
We will denote it by $s({\cal L})$. 
Let 
$L_3\subset L_4 \subset \cdots \subset L$ be a filtration, where 
$\dim L_j = j$ runs from $3$ to $\dim L$.  
The invariant $\dim \Sigma_\eta$  
defines a decomposition of $\Lambda^3L/GL(L)$ as a union 
$$\{0 \}\cup \bigcup_{j=3}^{\dim L}\Lambda^3_{reg}L_j/GL(L_j),$$ 
where $\Lambda^3_{reg}\bR^j \subset \Lambda^3\bR^j$ is the open subset
of 3-vectors with $j$-dimensional support. The points of the stratum 
$\Lambda^3_{reg}L_j/GL(L_j)\cong \Lambda^3_{reg}\bR^j/GL(j)$ correspond
to isomorphism classes of pairs $({\cal L},h)$ with $s({\cal L})=j$.\\ 
2) Since in the above classification $\Sigma_\eta$ is isotropic, 
it is clear that a flat (or 2-step nilpotent) bi-invariant metric on a  
Lie group is indefinite, unless $\eta=0$ and the group is Abelian. 
It follows from Milnor's classification of Lie
groups with a flat left-invariant Riemannian metric \cite{Mi} that
a 2-step nilpotent Lie group with a flat left-invariant Riemannian metric
is necessarily Abelian.\\

Since a nilpotent Lie group with rational structure constants has a 
(co-compact) lattice \cite{Ma}, we obtain: 

\bc \label{Cor} The groups $({\cal L}(\eta ),h)$ admit 
lattices $\Gamma \subset {\cal L}(\eta )$,  
provided that $\eta$ has rational coefficients with 
respect to some basis.  
$M=M(\eta, \Gamma ):= \Gamma \setminus {\cal L}(\eta )$ is a 
{\em flat compact homogeneous 
pseudo-Riemannian manifold}. The connected 
component of the identity in the isometry group of $M$ 
is the image of the natural group homomorphism $\pi$ from 
${\cal L}(\eta )$ into the isometry group of 
$M$. 
\ec  

\pf 
First we remark that the  bi-invariant metric $h$ induces an 
${\cal L}(\eta )$-invariant metric on the homogeneous space 
$M= \Gamma \setminus {\cal L}(\eta )$. We shall identify  
the group $\Gamma$ with a subgroup of the  isometry group of
$\widetilde{M}:=({\cal L}(\eta ),h)$ using the action of $\Gamma$ on 
${\cal L}(\eta )$ by left-multiplication. 
Let $G$ be the connected component of the identity in the isometry group of
$\widetilde{M}$.  It is clear that any element of $G$ which commutes 
with the action of $\G$ induces an isometry of $M$. Therefore  
we have a natural homomorphism $Z_G(\Gamma ) \ra \mathrm{Isom}(M)$ from 
the centraliser $Z_G(\Gamma )$  of $\Gamma$ in $G$ into $\mathrm{Isom}(M)$.  
In particular, the connected group $Z_G(\Gamma )_0$ is mapped into 
$\mathrm{Isom}_0(M)$.
Conversely, the action of $\mathrm{Isom}_0(M)$ on $M$ can be lifted 
to the action of a connected Lie 
subgroup $H\subset G$ on $\widetilde{M}$,
which maps cosets of $\G$ to cosets of $\G$. The latter property 
implies that $H$ normalises the subgroup $\G \subset 
\mathrm{Isom}(\widetilde{M})$.   
Since $\G$ is discrete and $G$ is connected, we can conclude that 
$H$ is a subgroup of the centraliser $Z_G(\Gamma )$  of $\Gamma$ in $G$. 
As $H$ is connected, we obtain $H\subset Z_G(\Gamma )_0$. 
By the previous argument, we have also the opposite inclusion 
$Z_G(\Gamma )_0 \subset H$ and, hence, 
$H=Z_G(\Gamma )_0$. 
Now the statement about the isometry group of $M$ follows from the fact that 
the centraliser in $G$ 
of the left-action of $\Gamma \subset 
{\cal L}(\eta )$ on ${\cal L}(\eta )$ 
is precisely the right-action of ${\cal L}(\eta )$ on ${\cal L}(\eta )$,
since $\Gamma \subset {\cal L}(\eta )$ is Zariski-dense, see \cite{R} 
Theorem 2.1. In fact, this shows that $H$ coincides with the 
group ${\cal L}(\eta )$
acting by right-multiplication on $\widetilde{M}={\cal L}(\eta )$
and that $\mathrm{Isom}_0(M)$ coincides with ${\cal L}(\eta )$
acting by right-multiplication on $M=\Gamma \setminus {\cal L}(\eta )$. 
\qed 

\begin{Ex}
We consider $V=(\bR^{3,3}, \langle \cdot , \cdot \rangle) $ and a basis
$(e_1,e_2,e_3,f_1,f_2,f_3 )$  
such that $ \langle e_i , f_j \rangle= \delta_{ij}$ and  $ \langle e_i , e_j \rangle=  \langle f_i , f_j \rangle=0.$
Then the three-vector $\eta := f_1\wedge f_2\wedge f_3\in \wedge^3V$ has 
isotropic support $\Sigma_{\eta}=\mbox{span}\{f_1,f_2,f_3 \}$. 
The non-vanishing components of the Lie bracket defined 
by \re{LiebracketEqu} are  
$$[e_1,e_2]=2 f_3, [e_2,e_3]=2f_1, [e_3,e_1]=2 f_2.$$ 
We have seen above that the bi-invariant metric $h$ was obtained by
pulling back the scalar product 
$\langle \cdot , \cdot \rangle$ by the orbit map
\eqref{def_obitmap} which identifies $\mathcal{L}(\eta)$ with $V$ 
via $ \mathcal{L}(\eta) \ni g_X \mapsto g_X0=X \in V.$ The inverse map is 
$V \ni X \mapsto g_X \in  \mathcal{L}(\eta).$ This identifies 
the pseudo-Riemannian 
manifolds $(\mathcal{L}(\eta),h)$ and $(V, \langle \cdot , \cdot \rangle).$ 
In consequence the isometry group
of $\mathcal{L}(\eta)$ is isomorphic to the full affine pseudo-orthogonal
group operating by $g_X \mapsto g_{AX+ v}$ with $A \in O(V)$ and $v \in V.$ 
Next we consider the lattice 
\[ \G := \{ g_Y|Y\in \bZ^6\}, \]
where $\bZ^6 \subset V$ is the lattice of integral vectors with 
respect to the basis  $(e_1,e_2,e_3,f_1,f_2,f_3 )$. 
An element 
$g_Y \in \G$ operates from the left on $\mathcal{L}(\eta)\cong V$ as 
$$ X \mapsto (\id + \eta_{Y})X +  Y. $$ 
Let us determine the centraliser of this $\G$-action in the
isometry group of $\mathcal{L}(\eta).$ A short calculation shows that 
an affine isometry $(A,v)$ with linear part 
$A \in O(V)$ and translational part  $v\in V$ belongs to the centraliser 
of $\G$ if and only if 
\[ [\eta_Y,A]X+\eta_Yv - AY + Y =0\]
for all $X\in V$, $Y\in \bZ^6$. For $X=0$ we get 
$AY=\eta_Yv+Y=(\id - \eta_v)Y$ and, hence, $A=\id - \eta_v$. This 
shows that the affine transformation $(A,v)$ corresponds to 
the right action of the element $g_v$, which obviouly belongs to the
centraliser. Therefore, in this example, 
we have proven by direct calculation that the
centraliser in the isometry group of $\mathcal{L}(\eta)$ 
of $\G$ acting by left-multiplication on $\mathcal{L}(\eta)$  
is precisely the group $\mathcal{L}(\eta)$ acting by right-multiplication. 
This fact was proven for arbitrary groups $\mathcal{L}(\eta)$ 
and lattices $\G$ in the proof of Corollary \ref{Cor}. 
\end{Ex}
 
\section{Flat nearly para-K\"ahler manifolds}
In this section we give a constructive classification of 
flat nearly para-K\"ahler manifolds and show that such manifolds
provide a class of examples for the flat Lie groups discussed in section 
\ref{simplytransSec}. The structure of the section is as follows. In the
first subsection we give a short introduction to para-complex geometry. For
more information the reader is referred to \cite{CMMS}. The second part 
discusses nearly para-K\"ahler manifolds and derives some consequences of 
the flatness. In the third subsection we give a local classification which
relates a flat nearly para-K\"ahler manifold to an element of a certain subset
${\cal C}_\tau(V)$ of the cone ${\cal C}(V)\subset \wedge^3V$ defined
in Theorem \ref{2stepnp_Thm}.  The structure of ${\cal C}_\tau(V)$ is studied 
in the last subsection and global classification results are derived.

\subsection{Para-complex geometry}
The idea  of para-complex geometry is to replace the complex structure $J$ 
satisfying $J^2=- Id$  on a (finite) dimensional vector space $V$ by a {\cmssl para-complex structure} $\I$ satisfying $\I^2 =  Id$  
and to require that the two eigenspaces of $\I,$ i.e. $V^{\pm}:= \ker(Id \mp \I),$ have the same dimension. A 
{\cmssl para-complex vector space} $(V,\I)$ is a vector space endowed with a para-complex structure. 
Para-complex, para-Hermitian and para-K\"ahler geometry was first studied in \cite{L}. We
invite the reader to consult  \cite{CFG} or the more recent article \cite{AMT}
for a survey on this subject.
\bd
An {\cmssl almost para-complex structure} $\I$ on a smooth manifold $M$ is an 
endomorphism field $\I \in \G(\mbox{End}(TM), p \mapsto \I_p,$ such that 
$\I_p$ is a para-complex structure on $T_pM$ for all points $p\in M.$ A manifold endowed 
with an almost para-complex structure is called an {\cmssl almost para-complex manifold}.  \\
An almost para-complex structure is called {\cmssl integrable} 
if its eigendistributions
$T^{\pm}M:= \ker(Id \mp \I)$ are both integrable. A manifold endowed 
with an integrable almost para-complex structure is called a {\cmssl  para-complex manifold}.
\ed
We remark that the obstruction to integrability (cf.\ Proposition 1 of \cite{CMMS}) of an almost para-complex structure is the {\cmssl Nijenhuis tensor} of $\tau$, which is the tensor field defined by 
$$ N_\I (X,Y):= [X,Y] +[\I X,\I Y]- \I[X,\I Y] - \I[\I X,Y],$$
for all vector fields $X, Y$ on $M$. 
\bd
Let $(V,\I)$ be a para-complex vector space. A {\cmssl para-Hermitian scalar 
product} $g$ 
on $(V,\I)$ is a pseudo-Euclidian scalar product, such that 
$\I^*g(\cdot,\cdot)= g(\I\cdot,\I\cdot)= -g(\cdot,\cdot).$\\
A {\cmssl para-Hermitian vector space} is a para-complex vector space endowed with a 
para-Hermitian scalar product. The pair $(\I,g)$ is called 
{\cmssl para-Hermitian structure} on the vector space $V.$
\ed
The next two examples give two frequently used models of para-Hermitian structures:
\begin{Ex}
Let us consider the vector space $\bR^{2n}=\bR^{n}\oplus\bR^{n}$
and denote by $e^+_i=e_i\oplus 0$ and $e^-_i=0\oplus e_i$ its standard basis.
Its standard para-complex structure is given by $\I e^\pm_i=\pm e^\pm_i.$ A para-Hermitian 
scalar product $g$ is given by $g(e^\pm_i,e^\pm_j)=0$ and $g(e^\pm_i,e^\mp_j)=\delta_{ij}.$ 
We call the pair $(\I,g)$ the {\cmssl standard para-Hermitian structure} of $\bR^{2n}.$
\end{Ex}
\begin{Ex}
We denote by $C=\bR [e]\cong \bR \oplus \bR$, $e^2=1$, the ring of 
para-complex numbers. 
Consider the real vector space $C^n=\bR^n\oplus \, e \bR^n$ with  
standard basis given by 
$(e_1,\ldots,e_n,f_1,\ldots,f_n),$ where $f_i=ee_i$ and its standard para-complex structure which 
is defined by $\I e_i = f_i$ and $\I f_i = e_i.$ Then we can define a para-Hermitian scalar 
product by $g(e_i,e_j)=-g(f_i,f_j)=\delta_{ij}$ and $g(e_i,f_j)=0.$ We 
denote this pair $(\I,g)$ the  {\cmssl standard para-Hermitian structure} of $C^{n}.$ 
\end{Ex}
The decomposition of the cotangent bundle $T^*M= (T^*M)^+ \oplus (T^*M)^-$
with respect to the dual para-complex structure induces a bi-grading on the bundle of 
exterior forms $ \Lambda^k T^*M = \oplus_{k=p+q}\, \Lambda^{p,q}T^*M.$
An element of $\Lambda^{p,q}T^*M$ will be called of {\cmssl type} $(p,q).$ The corresponding 
decomposition on differential forms is denoted by $ \O^k(M)= \oplus_{k=p+q} \,\Omega^{p,q}(M).$
\bd
An {\cmssl almost para-Hermitian manifold} $(M,\I ,g)$ 
is an almost para-complex 
manifold $(M,\I)$ which is endowed with a pseudo-Riemannian metric $g$ 
which is {\cmssl para-Hermitian}, i.e.\ it satisfies $\I^*g(\cdot,\cdot) = g(\I\cdot,\I \cdot)=-g(\cdot,\cdot).$ 
\ed 
Note that the condition on the metric to be para-Hermitian forces it to have split signature $(n,n).$

\subsection{Basic facts and results about nearly para-K\"ahler manifolds}  
The notion of a {\it nearly para-K\"ahler manifold} was recently introduced by Ivanov and Zamkovoy \cite{I}.
\bd 
An almost para-Hermitian manifold $(M,\I,g)$ is called {\cmssl nearly
  para-K\"ahler} manifold, if its Levi-Civita  connection $D$ satisfies the
equation \be (D_X\I)Y = -(D_Y\I)X, \quad \forall   X,Y \in
\Gamma(TM). \label{NK_def} \ee
A nearly para-K\"ahler manifold is called {\cmssl strict}, if $D \I \ne 0.$
\ed
Like for a nearly K\"ahler manifold there exists a canonical para-hermitian
connection with totally skew-symmetric torsion.
\bp[Prop. 5.1 in \cite{I}] \label{Can_con} Let $(M,\I,g)$ be a nearly para-K\"ahler manifold. 
Then  there exists a unique connection $\n$ with totally skew-symmetric 
torsion $T^{\n}$ (i.e. $g(T^{\n}(\cdot,\cdot),\cdot)$ is a three-form)  satisfying $\n g=0$ and $\n \I =0.$  \ep
\noindent
More precisely, this connection is given by
\be \n_XY  = D_XY - \eta_XY\mbox{ with } \eta_XY = -\frac{1}{2} \I (D_X \I)Y
\mbox{ and } X,Y \in \Gamma(TM) \ee
and in consequence the torsion is
\be T^{\nabla} = -2 \eta  \ee
and one has $\{\eta_X,\I \} =0$ for all vector fields $X.$
In the same reference \cite{I} Theorem 5.3 it is shown that, as in the nearly
K\"ahler case, the torsion of $\n$ is parallel, i.e.
 \be \label{tor_n_para} \n \eta=0 \mbox{ and } \n (D\I)=0. \ee
\bp \label{eta_D_nabla_prop}
Let $(M,g,\I)$ be a flat nearly para-K\"ahler manifold, then
\begin{enumerate}
\item[1)] $ \eta_X \circ \eta_Y =0$ for all $X,Y,$
\item[2)]$D \eta =\n \eta =0.$ 
\end{enumerate}
\ep
\pf
On a nearly para-K\"ahler manifold one has the identity $$R^D(X,Y,Z,W)
+R^D(X,Y,\I Z,\I W) = g((D_X\I))Y,(D_Z\I)W),$$ cf. \cite{I} Proposition
5.2. For a flat nearly para-K\"ahler manifold it follows
\bea
 g((D_X\I)Y,(D_Z\I)W)=0 \quad \forall X,Y,Z,W.
\eea
\noindent
With this identity  and  $D\I \circ \I = -\I \circ D \I$ we obtain
\bean
0=g((D_X\I )Y,(D_Z\I )W)&=&-g((D_Z\I )(D_X\I )Y,W)
= 4g(\eta_Z\circ \eta_X Y,W).
\eean
This shows $\eta_X \circ \eta_Y =0$ for all $X,Y$ and finishes the proof of 
part 1.). \\
2.) With two vector fields $X,Y$ we calculate
\bean
(D_X\eta)_Y &=& D_X(\eta_Y)- \eta_{D_XY} \overset{D= \n +\eta}{=}\n_X(\eta_Y) + [\eta_X ,\eta_Y] - \eta_{D_XY}\\ 
&=& (\n_X \eta)_Y + \eta_{[\n_XY -D_XY]} + [\eta_X,\eta_Y]= (\n_X \eta)_Y - \eta_{\eta_XY } + [\eta_X,\eta_Y]\\
&\overset{\eqref{NK_def}}{=}& (\n_X \eta)_Y + \eta_{\cdot} \eta_XY  + [\eta_X,\eta_Y] \overset{1.)}{=} (\n_X \eta)_Y\overset{\eqref{tor_n_para}}=0.
\eean This is part 2).
\qed

\subsection{Local classification of flat nearly para-K\"ahler manifolds}
We consider $(C^n,\I_{can})$ endowed with the standard $\I_{can}$-anti-invariant pseudo-Euclidian scalar product $g_{can}$ of signature $(n,n).$

Let $(M,g,\I )$ be a flat nearly para-K\"ahler manifold. Then there exists for each point $p \in M$ an open set
 $U_p \subset M$ containing the point $p,$ a connected open set $U_0$ of $C^n$ containing the origin $0 \in C^n$ and an isometry 
$ \Phi\,:\, (U_p,g) \tilde{\rightarrow} (U_0,g_{can}),$  such that in $p  \in M$ we have $\Phi_* \I_p = \I_{can} \Phi_*.$
In other words, we can suppose, that locally $M$ is a connected  open subset
of $C^n$ containing the origin $0$ and that $g=g_{can}$ and $\I_0=\I_{can}.$
Summarizing Proposition \ref{Can_con} and \ref{eta_D_nabla_prop} we obtain the
next Corollary.
\bc \label{class_cor}
Let $M \subset  C^n$ be an open neighborhood of the origin endowed with a nearly para-K\"ahler structure $(g,\I)$ such 
that $g=g_{can}$ and $\I_0=\I_{can}$. The $(1,2)$-tensor $$ \eta := -\frac{1}{2} \I  D\I  $$ defines a constant three-form on $M \subset C^n =\bR^{n,n}$ given 
by $\eta(X,Y,Z)= g(\eta_XY,Z) $ and satisfying 
\begin{enumerate}
\item[(i)]  $\eta \in {\cal C}(V),$ i.e. $\eta_X \, \eta_Y =0,\quad \forall X,Y,$
\item[(ii)] $\{\eta_X,\I_{can} \}=0, \quad \forall X.$ 
\end{enumerate}
\ec
The rest of this subsection is devoted to the local classification
result. In subsection \ref{glob_class} we study the
structure of the subset of ${\cal C}(V)$ given by the condition (ii) in more
detail and give global classification results.  The converse statement of Corollary \ref{class_cor} is given in the next lemma.
\bl \label{Lemma_J_from_eta}
Let $\eta$ be a constant three-form on an open connected set  $M \subset C^n$ of $0$ satisfying (i) and (ii) of Corollary \ref{class_cor}. Then there 
exists a  unique para-complex structure $\I $ on $M$ such that
\begin{enumerate}
\item[a)] $\I_0=\I_{can},$
\item[b)] $\{\eta_X,\I\}=0,\quad \forall X,$
\item[c)] $D\I  =-2\I  \eta,$
\end{enumerate}
where $D$ is the Levi-Civita connection of the pseudo-Euclidian vector space $C^n.$ \\
Let $\n:=D-\eta$ and assume b) then c) is equivalent to 
\begin{enumerate}
\item[c)'] $ \n \I  =0.$ 
\end{enumerate}
Furthermore, this para-complex structure $\tau$ is skew-symmetric with respect to $g_{can}.$
\el
\pf
One proves the equivalence of c) and c)' by an easy computation. \\
Let us show the uniqueness: Given two almost para-complex structures
satisfying a)-c) we deduce $(\I -\I')_0=0$ and  $\n \I  =\n \I' =0.$ This shows $\I  \equiv \I'.$ To show the existence we define
\bea
 \I &=& \exp\left( 2 \sum_{i=1}^{2n} x^i\, \eta_{\partial_i}\right) \I_{can} \overset{(i)}{=} \left( Id + 2 \sum_{i=1}^{2n} x^i \, \eta_{\partial_i}\right)\I_{can}, \label{def_I}
\eea
where $x^i$ are linear coordinates of $C^n =\bR^{n,n} =\bR^{2n}$ and
$\partial_i =\frac{\partial}{\partial x^i}.$ \\
{\bf Claim:} $\I$ defines a para-complex structure which satisfies a)-c). \\
a) From $x^i(0)= 0$ we obtain $\I_0 =\I_{can}.$\\
b) Follows from the definition of $\I$ (cf. equation \eqref{def_I}) and the properties (i) and (ii).  \\
c) One computes 
\bean
D_{\partial_j} \I &=& 2 \exp \left( 2 \sum_{i=1}^{2n} x^i \, \eta_{\partial_i} \right) \eta_{\partial_j} \I_{can}\overset{(ii)}{=} -2 \underbrace{\exp \left( 2 \sum_{i=1}^{2n} x^i \, \eta_{\partial_i} \right)\,  \I_{can}}_\I\, \eta_{\partial_j}=-2 \I\, \eta_{\partial_j}.
\eean
It holds $ \I= \I_{can} + \left( 2 \sum_{i=1}^{2n} x^i\, \eta_{\partial_i}\right) \I_{can},$
where $\{ \eta_{\partial_i},\I_{can} \} =0$ and $ \eta_{\partial_i}$ is $g$-skew-symmetric. This implies that $\I$ is $g$-skew-symmetric. It remains to prove $\I^2=Id.$ 
\bean 
\I^2 &= & \left( Id + 2 \sum_{i=1}^{2n} x^i \, \eta_{\partial_i}\right)\I_{can} \left( Id + 2 \sum_{i=1}^{2n} x^i \, \eta_{\partial_i}\right)\I_{can}\\&=& \left( Id + 2 \sum_{i=1}^{2n} x^i \, \eta_{\partial_i}\right) \left( Id - 2 \sum_{i=1}^{2n} x^i \, \eta_{\partial_i}\right) = \left[ Id -4 \left(\sum_{i=1}^{2n} x^i \, \eta_{\partial_i}\right)^2\right] \overset{(i)}{= }Id.
\eean 
This finishes the proof of the lemma.
\qed
\bt \label{1stThm}
Let $\eta$ be a constant three-form on a connected open set $U \subset C^n$ containing the origin $0$ 
which satisfies  (i) and (ii) of Corollary \ref{class_cor}. Then there exists a unique almost para-complex structure  
\be
\I=  \exp\left( 2 \sum_{i=1}^{2n} x^i\, \eta_{\partial_i}\right) \I_{can}
\ee
on $U$ such that a) $\I_0 =\I_{can},$ and b) $M(U,\eta):= (U,g=g_{can},\I)$ is
a flat nearly para-K\"ahler manifold. Any flat nearly para-K\"ahler manifold is locally isomorphic to a flat nearly 
para-K\"ahler manifold of the form $M(U,\eta).$
\et
\pf
$(M,g)$ is a flat pseudo-Riemannian manifold. Due to Lemma  \ref{Lemma_J_from_eta} $\I,$ 
is a skew-symmetric almost para-complex structure on $M$ and $\I_0=\I_{can}.$ From Lemma \ref{Lemma_J_from_eta}  c) and the skew-symmetry of $\eta$ it follows the skew-symmetry of $D\I.$ Therefore $(M,g,\I)$ is a nearly para-K\"ahler manifold. 
The remaining statement follows from Corollary \ref{class_cor} and Lemma \ref{Lemma_J_from_eta}.
\qed 

\subsection{The variety $\mathcal{C}_{\I}(V)$} \label{glob_class}
Now we discuss the solution of (i) and (ii) of Corollary \ref{class_cor}. 
In the following we shall freely 
identify the real vector space $V:=C^n=\bR^{n,n}=\bR^{2n}$ with
its dual $V^*$ by means of the pseudo-Euclidian scalar product $g=g_{can}$. 
The geometric interpretation is given in terms of an affine variety $\mathcal{C}_{\I}(V)\subset \Lambda^3 V.$ 

\bp \label{char_i_prop}
A three-form $\eta \in \Lambda^3 V^*\cong \Lambda^3 V$ 
satisfies (i) of Corollary \ref{class_cor}, i.e. $\eta_X \circ \eta_Y=0, X,Y, \in V,$ if and only if 
there exists  an isotropic subspace $L\subset V$
such that $\eta \in \Lambda^3 L \subset \Lambda^3 V$. 
If $\eta$ satisfies (i) and (ii) of Corollary \ref{class_cor} then
there exists a $\I_{can}$-invariant isotropic subspace $L\subset V$ 
with $\eta \in \Lambda^3 L$. 
\ep

\pf The proposition follows from the next lemma by taking $L=\Sigma_\eta$. 
\qed 

\bl \label{Lemma_cond_I}\vspace{-1.3em}
\begin{enumerate}
\item
$\Sigma_{\eta}$ is isotropic if and only if $\eta$ satisfies (i) of  
Corollary \ref{class_cor}. If $\eta$ satisfies (ii) of  Corollary 
\ref{class_cor}, then
$\Sigma_{\eta}$ is $\I_{can}$-invariant. 
\item
 Let $\eta \in \Lambda^3V$.  Then $\eta \in \Lambda^3\Sigma_{\eta}.$ 
\end{enumerate}
\el 
\pf
The proof of the first part is analogous to Lemma 6 in  \cite{CS}. The 
second part is Lemma 7 of \cite{CS}.
\qed\\\noindent 
Any three-form $\eta$ on $(V,\I_{can})$ decomposes with respect to the 
grading induced by the decomposition $V = V^{1,0} \oplus V^{0,1} $
into $ \eta = \eta^+ + \eta^- $ with  $\eta^+ \in \Lambda^+V := \Lambda^{2,1}V + \Lambda^{1,2}V$ and  $\eta^- \in  \Lambda^-V := \Lambda^{3,0}V + \Lambda^{0,3}V.$
\bt \label{2ndThm} 
A three-form $\eta \in \Lambda^3 V^*\cong \Lambda^3 V$ satisfies (i)
and (ii) of Corollary \ref{class_cor} if and only if there exists  an 
isotropic $\I_{can}$-invariant subspace $L$ 
such that  $\eta \in \Lambda^-L=\Lambda^{3,0}L + \Lambda^{0,3}L \subset \Lambda^3L\subset \Lambda^3 V$ 
(The smallest such subspace $L$ is $\Sigma_\eta$.).  
\et  \vspace{-0.5em} \noindent
We need the following general Lemma. \vspace{-0.5em}
\bl 
It is 
\bean   \Lambda^-V = 
\{ \eta \in \Lambda^3 V\,|\, \eta(\cdot,\I\cdot,\I \cdot) = \eta(\cdot,\cdot,\cdot) \} 
= \{ \eta \in \Lambda^3 V\,|\,\{\eta_X,\I \}=0, \; \forall X \in V \}.
\eean
\el

\pf (of Theorem \ref{2ndThm})  By Proposition \ref{char_i_prop}, the conditions (i)
and (ii) of Corollary \ref{class_cor} imply the existence
of an isotropic $\I_{can}$-invariant subspace $L\subset V$ such that 
$\eta \in \Lambda^3L$. The last lemma shows that the condition (ii) 
is equivalent to $\eta \in \Lambda^-V.$ Therefore 
$\eta \in \Lambda^3L\cap  \Lambda^-V = \Lambda^-L$. 
The converse statement follows from the same argument. \qed 

\bc \label{cor_Th1_th2}
\begin{itemize}
\item[(i)] The conical affine variety $ \mathcal{C}_\I(V):= \{ \eta\, | \,
  \eta \mbox{ satisfies } (i) \mbox{ and } (ii) \mbox{ in Corollary
    \ref{class_cor} }\}\subset \Lambda^3V$
has the following description
$ \mathcal{C}_\I(V)= \underset{L\subset V}{\bigcup} \Lambda^-L= \underset{L\subset V}{\bigcup} (\Lambda^3L^+ + \Lambda^3L^-), $
where the union is over all $\I$-invariant maximal isotropic subspaces.
\item[(ii)]
If $\dim V < 12$ then it holds 
$ \mathcal{C}_\I(V) = \Lambda^3V^+ \cup \Lambda^3V^-.$
\item[(iii)]
Any  flat nearly para-K\"ahler  manifold $M$ is locally of the form $M(U,\eta )$, for some $\eta \in \mathcal{C}_\I(V)$  
and some open subset $U \subset V.$
\item[(iv)]
There are no strict  flat nearly para-K\"ahler manifolds of 
dimension less than 6.
\end{itemize}
\ec
\pf (i) follows from Theorem \ref{2ndThm}. \\
(ii) Let $L \subset V$ be a $\I$-invariant isotropic subspace. If $\dim V < 12,$ then
$\dim L < 6$ and, hence, either $\dim L^+ < 3$ or $\dim L^- < 3.$ In the first case we have
$$ \Lambda^-L = \Lambda^3L^+ + \Lambda^3L^- =\Lambda^3L^- \subset \Lambda^3V^-,$$ 
in the second case it is $ \Lambda^-L = \Lambda^3L^+ + \Lambda^3L^- =\Lambda^3L^+ \subset \Lambda^3V^+.$\\
(iii)  is a consequence of (i), Theorem \ref{1stThm} and \ref{2ndThm}. \\
(iv) By (iii) the strict nearly para-K\"ahler manifold $M$ is locally of the form $M(U,\eta ),$ which is strict if and only 
if $\eta\neq 0$. This is only possible for 
$\dim L \ge 3$, i.e.\ for $\dim M\ge 6$. \qed

\begin{Ex}
We have the following example which shows that part (ii) of 
Corollary \ref{cor_Th1_th2} fail in dimension $\ge$ 12: Consider $(V,\tau)=(C^6,e)=\bR^6\oplus \, e \bR^6$ with a basis given by 
$(e^+_1,\ldots,e^+_6,e^-_1,\ldots,e^-_6),$ such that $e^\pm_i$ form a 
basis of $V^\pm$ with $g(e_i^+,e_j^-)= \delta_{ij}.$ Then the form $\eta := e^+_1\wedge e^+_2\wedge e^+_3 + e^-_4\wedge e^-_5\wedge e^-_6$
lies in the variety $ \mathcal{C}_\I(V).$ 
\end{Ex}

\bt \label{3rdThm} Any strict flat nearly para-K\"ahler manifold 
is locally a pseudo-Riemannian product $M=M_0\times M(U,\eta )$ of a flat 
para-K\"ahler factor $M_0$ 
of maximal dimension and a flat nearly para-K\"ahler  
manifold $M(U,\eta )$, $\eta\in C_\tau (V)$, of signature $(m,m)$, 
$2m=\dim M(U,\eta )\ge 6$ such that $\Sigma_\eta$ has dimension  $m$. 
\et 

\pf  By Theorem \ref{1stThm} and \ref{2ndThm}, $M$ is locally isomorphic 
to an open subset of a manifold of 
the form $M(V,\eta )$, where $\eta \in \Lambda^3V$ has a $\I_{can}$-invariant 
and isotropic support $L=\Sigma_\eta$. We choose a   $\I_{can}$-invariant
isotropic subspace $L'\subset V$ such that $V':=L+L'$ is nondegenerate and 
$L\cap L'=0$ and put $V_0 = (L+L')^\perp$. Then $\eta \in \Lambda^3V'
\subset  \Lambda^3V$  and $M(V,\eta ) = M(V_0,0)\times M(V',\eta )$. 
Notice that $M(V_0,0)$ is simply the flat para-K\"ahler manifold $V_0$ 
and that  $M(V',\eta )$ is strict 
of split signature $(m,m)$, where 
$m=\dim L\ge 3$. \qed \\

\bc \label{lastCor} Any simply connected nearly para-K\"ahler manifold with a 
(geodesically) complete flat metric  
is a pseudo-Riemannian product $M=M_0\times M(\eta )$ of a flat 
para-K\"ahler factor $M_0=\bR^{l,l}$ 
of maximal dimension and a flat nearly para-K\"ahler  
manifold $M(\eta ):=M(V,\eta )$, $\eta\in C_\tau (V)$, of signature $(m,m)$ 
such that $\Sigma_\eta$ has dimension  
$m=0,3,4,\ldots$. 
\ec 

\noindent
Next we wish to describe the moduli space of (complete simply connected) flat  
nearly para-K\"ahler manifolds $M$ of dimension $2n$ up to isomorphism. 
Without
restriction of generality we will assume that $M=M(\eta )$ 
has no para-K\"ahler 
de Rham factor, which means that $\eta\in C_\tau (V)$ has maximal 
support $\Sigma_\eta$, i.e.\ $\dim \Sigma_\eta = n$. 
We denote by $C_\tau^{reg}(V)\subset C_\tau(V)$ the open 
subset consisting of elements with maximal support. 
The group 
\[ G:={\rm Aut}(V,g_{can},\tau_{can}) \cong GL(n)\]
acts on  $C_\tau(V)$ and preserves $C_\tau^{reg}(V)$.  
Two nearly para-K\"ahler manifolds $M(\eta )$ and $M(\eta')$ are
isomorphic if and only if $\eta$ and $\eta'$ are related by an element
of the group $G$. 
  
For $\eta\in C_\tau(V)$ we denote by $p$, $q$  the dimensions of
the eigenspaces of $\tau$ on $\Sigma_\eta$ for the eigenvalues $1, -1$, 
respectively. We call the pair $(p,q)\in \bN_0\times \bN_0$
the  {\cmssl type} of $\eta$. We will also say that the 
corresponding flat nearly para-K\"ahler manifold $M(\eta )$ has type $(p,q)$. 
We denote by $C^{p,q}_\tau(V)$ the 
subset of $C_\tau (V)$ consisting of elements of type $(p,q)$. 
Notice that $p+q\le n$ with equality if and 
only if $\eta \in C_\tau^{reg}(V)$. We have the following
decomposition 
$$C_\tau^{reg}(V) = \underset{(p,q)\in \Pi}{\bigcup}C^{p,q}_\tau (V),$$
where $\Pi := \{ (p,q)| p,q \in \bN_0\setminus \{1,2\}, p+q=n\}$. 
The group $G=GL(n)$ acts on the 
subsets $C^{p,q}_\tau(V)$ and we are interested in the orbit space
$C^{p,q}_\tau(V)/G$.  

Fix a $\tau$-invariant maximally isotropic subspace $L\subset V$
of type $(p,q)$ and put $\Lambda^-_{reg}L:=\Lambda^-L\cap C_\tau^{reg}(V)
\subset C^{p,q}_\tau (V)$. The stabilizer $G_L\cong GL(L^+)\times GL(L^-)
\cong GL(p)\times GL(q)$ of 
$L=L^++L^-$ in $G$ acts on $\Lambda^-_{reg}L$.

\bt
There is a natural one-to-one correspondence between complete simply connected 
flat nearly para-K\"ahler manifolds of type $(p,q)$, $p+q=n$, 
and the points of the following orbit space:  
$$C^{p,q}_\tau (V)/G\cong  \Lambda^-_{reg}L/G_L \subset \Lambda^-L/G_L =
\Lambda^3L^+/GL(L^+) \times \Lambda^3L^-/GL(L^-).$$
\et 

\pf 
Consider two complete simply connected 
flat nearly para-K\"ahler manifolds $M$, $M'$. By the previous
results we can assume that $M=M(\eta )$, $M'=M(\eta')$ are  
associated with $\eta, \eta' \in C^{p,q}_\tau (V)$. It is clear that 
$M$ and $M'$ are isomorphic if $\eta$ and
$\eta'$ are related by an element of $G$. To prove the converse
we assume that $\varphi : M \ra M'$ is an isomorphism of
nearly para-K\"ahler manifolds. By the results of Section \ref{simplytransSec}
$\eta$ defines a simply transitive group  of isometries. This group
preserves also the para-complex structure $\tau$, which is 
$\nabla$-parallel and hence left-invariant. This shows that $M$ and $M'$ 
admit a transitive group of automorphisms. Therefore, we can assume that
$\varphi$ maps the origin in $M=V$ to the origin in $M'=V$.  
Now $\varphi$ is an isometry of pseudo-Euclidian vector spaces
preserving the origin. Thus $\varphi$ is an element of $O(V)$
preserving also the para-complex structure $\I$ and hence
$\varphi \in G$. 

The identification of orbit spaces can be easily checked using Lemma \ref{Lemma_cond_I} 2.\ and the fact 
that any $\tau$-invariant isotropic subspace $\Sigma = \Sigma^++\Sigma^-$ can be mapped onto  $L$ by an element of $G$. 
\qed

\end{document}